\documentclass[11pt, reqno]{amsart} \textwidth=14.9cm \oddsidemargin=1cm \evensidemargin=1cm
\usepackage{amsmath,amssymb,amscd,mathrsfs,stmaryrd,wasysym}
\usepackage{amsmath}
\usepackage{amsxtra}
\usepackage{amscd}
\usepackage{amsthm}
\usepackage{amsfonts}
\usepackage{amssymb}
\usepackage{eucal}
\usepackage{color}
\usepackage{graphicx}%
\usepackage{bm}
\newtheorem{theorem}{Theorem}[section]
\newtheorem{lemma}[theorem]{Lemma}
\newtheorem{proposition}[theorem]{Proposition}

\theoremstyle{definition}
\newtheorem{definition}[theorem]{Definition}
\newtheorem{remark}[theorem]{Remark}
\newtheorem{example}[theorem]{Example}

\definecolor{A}{rgb}{.75,1,.75}

\numberwithin{equation}{section}

\begin{document}

\title[degenerate cyclotomic Yokonuma-Hecke algebras]{A note on degenerate cyclotomic Yokonuma-Hecke algebras}
\author[Weideng Cui]{Weideng Cui}
\address{School of Mathematics, Shandong University, Jinan, Shandong 250100, P.R. China.}
\email{cwdeng@amss.ac.cn}

\begin{abstract}
In this note, we first prove that the degenerate cyclotomic Yokonuma-Hecke algebra $Y_{r,n}^{d}$ is a cellular algebra by constructing an explicit cellular basis. We then develop the fusion procedure for $Y_{r,n}^{d}$, that is, we prove that the primitive idempotents of $Y_{r,n}^{d}$ can be defined by consecutive evaluations of a certain rational function.
\end{abstract}



\maketitle
\medskip
\section{Introduction}
\subsection{}
Yokonuma \cite{Yo} defined Yokonuma-Hecke algebras as a centralizer algebra associated to the permutation representation of a finite Chevalley group $G$ with respect to a maximal unipotent subgroup of $G$. The Yokonuma-Hecke algebra $Y_{r,n}(q)$ of type $A$ is defined as a quotient of the group algebra of the modular framed braid group $(\mathbb{Z}/r\mathbb{Z})\wr B_{n},$ where $B_{n}$ is the braid group on $n$ strands of type $A$. By the presentation given by Juyumaya and Kannan \cite{Ju1, Ju2, JuK}, the Yokonuma-Hecke algebra $Y_{r,n}(q)$ can also be regraded as a deformation of the group algebra of the complex reflection group $G(r,1,n),$ which is isomorphic to the wreath product $(\mathbb{Z}/r\mathbb{Z})\wr \mathfrak{S}_{n}$, where $\mathfrak{S}_n$ is the symmetric group.

Recently, by generalizing the approach of Okounkov-Vershik \cite{OV} on the representation theory of $\mathfrak{S}_n$, Chlouveraki and Poulain d'Andecy \cite{ChPA1} defined the affine Yokonuma-Hecke algebra $\widehat{Y}_{r,n}(q)$ and gave explicit formulae for all irreducible representations of $Y_{r,n}(q)$ over $\mathbb{C}(q)$, and further obtained a semisimplicity criterion for $Y_{r,n}(q)$. In the subsequent paper [ChPA2], they studied the representation theory of $\widehat{Y}_{r,n}(q)$ and the cyclotomic Yokonuma-Hecke algebra $Y_{r,n}^{d}(q)$. In particular, they gave the classification of irreducible representations of $Y_{r,n}^{d}(q)$ in the generic semisimple case. In \cite{CW}, we gave the classification of the simple $\widehat{Y}_{r,n}(q)$-modules as well as the classification of the simple modules of the cyclotomic Yokonuma-Hecke algebras over an algebraically closed field $\mathbb{K}$ of characteristic $p$ when $p$ does not divide $r.$ We \cite{C1} prove that the cyclotomic Yokonuma-Hecke algebra $Y_{r,n}^{d}(q)$ is cellular by constructing an explicit cellular basis. In the past few years, the study of affine and cyclotomic Yokonuma-Hecke algebras has made substantial progress; see [ChPA1-2, ChS, C1-6, CW, ER, JaPA, Lu, PA2, Ro].

\subsection{}
Jucys [Juc] claimed that the primitive idempotents of $\mathfrak{S}_n$ can be obtained by taking a certain limiting process on a rational function, which is now commonly referred to as the fusion procedure. It has been developed in the situation of Hecke algebras [Ch]; see also [Na1-3]. Molev [Mo] has proposed an alternative approach of the fusion procedure for the symmetric group, which is based on the existence of a maximal commutative subalgebra generated by the Jucys-Murphy elements. Here the idempotents can be derived by consecutive evaluations of a certain rational function. The version of the fusion procedure has been generalized to various algebras and groups; see [IMO, IM, IMOg1-2, OgPA1-2, PA1, ZL, C5-6].

\subsection{}
The purposes of this note are two folds. Firstly, we prove that the degenerate cyclotomic Yokonuma-Hecke algebra $Y_{r,n}^{d}$ is cellular by constructing an explicit cellular basis following the approach in [C1]. Secondly, we develop the fusion procedure for $Y_{r,n}^{d}$, that is, we prove that a complete set of pairwise orthogonal primitive idempotents of $Y_{r,n}^{d}$ can be constructed by consecutive evaluations of a certain rational function.

This paper is organized as follows. In Section 2, we first recall some combinatorial notions and then present the construction of a cellular basis of
$Y_{r,n}^{d}.$ In Section 3, for the split semisimple $Y_{r,n}^{d}$, we first deduce some properties of the idempotents $E_{\mathfrak{t}}$ of $Y_{r,n}^{d}$ by applying the general theory developed in [Ma, Section 3]. Then we give the inductive formulae for the idempotents $E_{\mathfrak{t}}$ in terms of the Jucys-Murphy elements. Finally we prove the fusion formulae for $E_{\mathfrak{t}}.$

\section{Cellular bases}

In this section, we first recall some combinatorial notions and then construct a cellular basis of the degenerate cyclotomic Yokonuma-Hecke algebra $Y_{r,n}^{d}.$

\subsection{Degenerate cyclotomic Yokonuma-Hecke algebras}
Let $r\in \mathbb{Z}_{\geq 1}$ and $p\in \mathbb{Z}_{\geq 1}$ such that $p$ does not divide $r$. Let $\mathbb{K}$ be an algebraically closed field of characteristic $p$ which contains some elements $v_{1},\ldots, v_{d}.$ The degenerate affine Yokonuma-Hecke algebra, denoted by $\widehat{Y}_{r,n}$, is an associative $\mathbb{K}$-algebra generated by the elements $t_{1},\ldots,t_{n},g_{1},\ldots,g_{n-1}, x_1,\ldots, x_{n}$ in which the generators $t_{1},\ldots,t_{n},g_{1},$ $\ldots,g_{n-1}$ satisfy the following relations:
\begin{equation}\label{rel-def-Y1}\begin{array}{rclcl}
g_ig_j\hspace*{-7pt}&=&\hspace*{-7pt}g_jg_i && \mbox{for all $i,j=1,\ldots,n-1$ such that $\vert i-j\vert \geq 2$;}\\[0.1em]
g_ig_{i+1}g_i\hspace*{-7pt}&=&\hspace*{-7pt}g_{i+1}g_ig_{i+1} && \mbox{for all $i=1,\ldots,n-2$;}\\[0.1em]
t_it_j\hspace*{-7pt}&=&\hspace*{-7pt}t_jt_i &&  \mbox{for all $i,j=1,\ldots,n$;}\\[0.1em]
g_it_j\hspace*{-7pt}&=&\hspace*{-7pt}t_{s_i(j)}g_i && \mbox{for all $i=1,\ldots,n-1$ and $j=1,\ldots,n$;}\\[0.1em]
t_i^r\hspace*{-7pt}&=&\hspace*{-7pt}1 && \mbox{for all $i=1,\ldots,n$;}\\[0.2em]
g_{i}^{2}\hspace*{-7pt}&=&\hspace*{-7pt}1 && \mbox{for all $i=1,\ldots,n-1$,}
\end{array}
\end{equation}
together with the following relations concerning the generators $x_1,\ldots, x_{n}$:
\begin{equation}\label{rel-def-Y2}\begin{array}{rclcl}
x_{i}x_{j}\hspace*{-7pt}&=&\hspace*{-7pt}x_{j}x_{i};\\[0.1em]
g_{i}x_{i+1}\hspace*{-7pt}&=&\hspace*{-7pt}x_{i}g_{i}+e_{i};\\[0.1em]
g_{i}x_{j}\hspace*{-7pt}&=&\hspace*{-7pt}x_{j}g_{i} \qquad \qquad \quad\mbox{for all $j\neq i, i+1$;}\\[0.1em]
t_{j}x_{i}\hspace*{-7pt}&=&\hspace*{-7pt}x_{i}t_{j} \hspace{0.07cm}\qquad \qquad\quad\mbox{for all $i, j=1,\ldots,n$,}
\end{array}
\end{equation}
where $s_{i}$ is the transposition $(i,i+1)$ and for each $1\leq i\leq n-1$,
$$e_{i} :=\frac{1}{r}\sum\limits_{s=0}^{r-1}t_{i}^{s}t_{i+1}^{-s}.$$

\begin{remark}\label{rem-YH}
The degenerate affine Yokonuma-Hecke algebra $\widehat{Y}_{r,n}$ is in fact a special case of the wreath Hecke algebra $\mathcal{H}_{n}(G)$ defined in [WW, Definition 2.4] when $G=C_{r}$ is the cyclic group of order $r;$ see also [RS].
\end{remark}
\noindent By [WW, Theorem 2.8], $\widehat{Y}_{r,n}$ has a $\mathbb{K}$-basis
\begin{equation}\label{PBW-basis}
\big\{t_{1}^{\beta_1}\cdots t_{n}^{\beta_{n}}x_{1}^{\alpha_1}\cdots x_{n}^{\alpha_{n}}g_{w}\:|\:0\leq \beta_1,\ldots,\beta_{n}\leq r-1, \alpha_1,\ldots,\alpha_{n}\geq 0, w\in \mathfrak{S}_{n}\big\}.\end{equation}

We define the degenerate cyclotomic Yokonuma-Hecke algebra $Y_{r, n}^{d}$ to be the quotient:
\[Y_{r, n}^{d}=\widehat{Y}_{r,n}/\langle (x_1-v_1)\cdots (x_{1}-v_{d})\rangle.\]
The following proposition has been provided in [WW, Proposition 5.5], which can be proved by adapting the approach in [Kle, Section 7.5] as has been done for cyclotomic Yokonuma-Hecke algebras in [C1]. It has also been proved in [Ro, Proposition 5.2] by adapting the approach in [AK].
\begin{proposition}\label{CYPBW-basis}
$Y_{r, n}^{d}$ has a $\mathbb{K}$-basis
\begin{equation}
\big\{t_{1}^{j_1}\cdots t_{n}^{j_{n}}x_{1}^{i_1}\cdots x_{n}^{i_{n}}g_{w}\:|\:0\leq i_1,\ldots,i_{n}\leq d-1,0\leq j_1,\ldots,j_{n}\leq r-1, w\in \mathfrak{S}_{n}\big\}.\end{equation}
\end{proposition}

Let $\mathbf{s}=\{1,2,\ldots,n\}.$ Let $i, k\in \mathbf{s}$ and set $$e_{i,k} :=\frac{1}{r}\sum\limits_{s=0}^{r-1}t_{i}^{s}t_{k}^{-s}.$$ Note that $e_{i,i}=1,$ $e_{i,k}^{2}=e_{i,k}=e_{k,i},$ and that $e_{i,i+1}=e_{i}.$ For any nonempty subset $I\subseteq \mathbf{s}$ we define the following element $E_{I}$ by $$E_{I} :=\prod_{i, j\in I; i<j}e_{i,j},$$where by convention $E_{I}=1$ if $|I|=1.$

We also need a further generalization of this. We say that the set $A=\{I_{1}, I_{2},\ldots,I_{k}\}$ is a set partition of $\mathbf{s}$ if the $I_{j}$'s are nonempty and disjoint subsets of $\mathbf{s},$ and their union is $\mathbf{s}.$ We refer to them as the blocks of $A.$ We denote by $\mathcal{SP}_{n}$ the set of all set partitions of $\mathbf{s}.$ For $A=\{I_{1}, I_{2},\ldots,I_{k}\}\in \mathcal{SP}_{n}$ we then define $E_{A} :=\prod_{j}E_{I_{j}}.$

\subsection{$(r,d)$-partitions}
A finite sequence $\lambda=(\lambda_{1},\ldots,\lambda_{k})$ is called a partition of $n$ if $\lambda_{1}\geq \cdots\geq \lambda_{k}\geq 0$ and $\lambda_{1}+\cdots+\lambda_{k}=n.$ If $\lambda$ is a partition of $n,$ we write $\lambda\vdash n$ and set $|\lambda| :=n$. We associate a Young diagram to a partition $\lambda$, which is the set $$[\lambda] :=\{(i,j)\:|\:i\geq 1~\mathrm{and}~1\leq j\leq \lambda_{i}\}.$$ We shall regard $[\lambda]$ as a left-justified array of rows such that there exist $\lambda_{j}$ nodes in the $j$-th row for $j=1,\ldots,k.$ If a node $\theta$ is in row $x$ and column $y,$ we write $\theta=(x,y)$.

For $\lambda\vdash n,$ a node $\theta\in [\lambda]$ is called removable from $\lambda$ if we still obtain a Young diagram of a partition after removing $\theta$ from $[\lambda]$; a node $\theta'\notin [\lambda]$ is called addable to $\lambda$ if we still obtain a Young diagram of a partition after adding $\theta'$ to $[\lambda]$. The conjugate of a partition $\lambda$ is the partition $\lambda'=(\lambda'_1,\ldots,\lambda'_l)$, which is given by $$\lambda'_j :=\sharp\{i\:|\:1\leq i\leq k \text{ such that }\lambda_i\geq j\}.$$

Let $d\in \mathbb{Z}_{\geq 1}.$ A $d$-partition of $n$ is an ordered $d$-tuple $\bm{\lambda}=(\lambda^{(1)},\lambda^{(2)},\ldots,\lambda^{(d)})$ of partitions $\lambda^{(k)}$ such that $\sum_{k=1}^{d}|\lambda^{(k)}|=n.$ The combinatorial objects appearing in the representation theory of $\mathrm{Y}_{r,n}^{d}$ will be $(r,d)$-partitions. By definition, an $(r,d)$-partition of $n$ is an ordered $r$-tuple $\underline{\bm{\lambda}}=(\bm{\lambda}^{1},\ldots,\bm{\lambda}^{r})=((\lambda_{1}^{(1)},\ldots,\lambda_{d}^{(1)}),\ldots,(\lambda_{1}^{(r)},\ldots,\lambda_{d}^{(r)}))$ of $d$-partitions $(\lambda_{1}^{(k)},\ldots,\lambda_{d}^{(k)})$ ($1\leq k\leq r$) such that $\sum_{k=1}^{r}\sum_{j=1}^{d}|\lambda_{j}^{(k)}|=n.$ We denote by $\mathcal{P}_{r,n}^{d}$ the set of $(r,d)$-partitions of $n.$ We shall say that the $l$-th partition of the $k$-th $r$-tuple has position $(k,l).$

An $(r,d)$-node $\bm{\theta}=(\theta, k, l)$ consists of a node $\theta,$ an integer $k\in \{1,\ldots,r\},$ and an integer $l\in \{1,\ldots,d\}$. We shall call $k$ the $r$-position of $\bm{\theta}$, $l$ the $d$-position of $\bm{\theta}$ and $(k,l)$ the position of $\bm{\theta}$. For each $\underline{\bm{\lambda}}\in \mathcal{P}_{r,n}^{d},$ we shall denote by $[\underline{\bm{\lambda}}]$ the set of $(r,d)$-nodes such that the subset consisting of the $(r,d)$-nodes having position $(k,l)$ forms a usual Young diagram of the partition $\lambda_{l}^{(k)}$, for any $k\in \{1,\ldots,r\}$ and $l\in \{1,\ldots,d\}.$

For each $\underline{\bm{\lambda}}\in \mathcal{P}_{r,n}^{d},$ an $(r,d)$-node $\bm{\theta}=(\theta, k, l)\in [\underline{\bm{\lambda}}]$ is called removable from $\underline{\bm{\lambda}}$ if the node $\theta$ is removable from the partition of $\underline{\bm{\lambda}}$ with position $(k,l)$; an $(r,d)$-node $\bm{\theta}=(\theta, k, l)\notin [\underline{\bm{\lambda}}]$ is called addable to $\underline{\bm{\lambda}}$ if the node $\theta$ is addable to the partition of $\underline{\bm{\lambda}}$ with position $(k,l).$

For an $(r,d)$-node $\bm{\theta}=((a,b),k,l)$, we define $\text{cc}(\bm{\theta}) :=b-a,$ $\text{p}^{(r)}(\bm{\theta}) :=k$, $\text{p}^{(d)}(\bm{\theta}) :=l$ and the content $\text{c}(\bm{\theta}) :=b-a+v_{l}.$

\subsection{Standard $(r,d)$-tableaux}
Let $\underline{\bm{\lambda}}$ be an $(r,d)$-partition of $n.$ An $(r,d)$-tableau $\mathfrak{t}=((\mathfrak{t}_{1}^{(1)},\ldots,\mathfrak{t}_{d}^{(1)}),\ldots,(\mathfrak{t}_{1}^{(r)},\ldots,\mathfrak{t}_{d}^{(r)}))$ of shape $\underline{\bm{\lambda}}$ is obtained by placing each $(r,d)$-node of $[\underline{\bm{\lambda}}]$ by one of the integers $1,2,\ldots,n,$ allowing no repeats. We will call the number $n$ the size of $\mathfrak{t}$ and the $\mathfrak{t}_{l}^{(k)}$'s the components of $\mathfrak{t}.$ Each $(r,d)$-node $\bm{\theta}$ of $\mathfrak{t}$ is labelled by $((a, b), k, l)$ if it lies in row $a$ and column $b$ of the component $\mathfrak{t}_{l}^{(k)}$ of $\mathfrak{t}.$

An $(r,d)$-tableau of shape $\underline{\bm{\lambda}}$ is called row standard if the numbers increase along any row (from left to right) of each diagram in $[\underline{\bm{\lambda}}].$ An $(r,d)$-tableau of shape $\underline{\bm{\lambda}}$ is called standard if the numbers increase along any row (from left to right) and down any column (from top to bottom) of each diagram in $[\underline{\bm{\lambda}}].$ From now on, we denote by $\text{Std}(\underline{\bm{\lambda}})$ the set of all standard $(r,d)$-tableaux of size $n$ and of shape $\underline{\bm{\lambda}},$ which is endowed with an action of $\mathfrak{S}_{n}$ from the right by permuting the entries in each $(r,d)$-tableau.

For each $\underline{\bm{\lambda}}\in \mathcal{P}_{r,n}^{d},$ we denote by $\mathfrak{t}^{\underline{\bm{\lambda}}}$ the standard $(r,d)$-tableau of shape $\underline{\bm{\lambda}}$ in which $1,2,\ldots,n$ appear in increasing order from left to right along the rows of the first diagram, and then along the rows of the second diagram, and so on.

For any partition $\lambda=(\lambda_{1},\ldots,\lambda_{k})$ of $n,$ we define the Young subgroup $\mathfrak{S}_{\lambda} :=\mathfrak{S}_{\lambda_{1}}\times\cdots\times\mathfrak{S}_{\lambda_{k}}.$ For each $\underline{\bm{\lambda}}=((\lambda_{1}^{(1)},\ldots,\lambda_{d}^{(1)}),\ldots,(\lambda_{1}^{(r)},\ldots,\lambda_{d}^{(r)}))\in \mathcal{P}_{r,n}^{d},$ we have a Young subgroup
$$\mathfrak{S}_{\underline{\bm{\lambda}}} :=\mathfrak{S}_{\lambda_{1}^{(1)}}\times\cdots\times\mathfrak{S}_{\lambda_{d}^{(1)}}\times\cdots\times
\mathfrak{S}_{\lambda_{1}^{(r)}}\times\cdots\times\mathfrak{S}_{\lambda_{d}^{(r)}},$$
which is exactly the row stabilizer of $\mathfrak{t}^{\underline{\bm{\lambda}}}.$

For a row standard $(r,d)$-tableau $\mathfrak{s}$ of shape $\underline{\bm{\lambda}},$ let $d(\mathfrak{s})$ be the element of $\mathfrak{S}_{n}$ such that $\mathfrak{s}=\mathfrak{t}^{\underline{\bm{\lambda}}}d(\mathfrak{s}).$ Then $d(\mathfrak{s})$ is a distinguished right coset representative of $\mathfrak{S}_{\underline{\bm{\lambda}}}$ in $\mathfrak{S}_{n},$ that is, $l(wd(\mathfrak{s}))=l(w)+l(d(\mathfrak{s}))$ for any $w\in \mathfrak{S}_{\underline{\bm{\lambda}}}.$ In this way, we obtain a correspondence between the set of row standard $(r,d)$-tableaux of shape $\underline{\bm{\lambda}}$ and the set of distinguished right coset representatives of $\mathfrak{S}_{\underline{\bm{\lambda}}}$ in $\mathfrak{S}_{n}.$

We define a partial order on the set of $(r,d)$-partitions and standard $(r,d)$-tableaux.

\begin{definition}\label{Def-def}
Let $\underline{\bm{\lambda}}=((\lambda_{1}^{(1)},\ldots,\lambda_{d}^{(1)}),\ldots,(\lambda_{1}^{(r)},\ldots,\lambda_{d}^{(r)}))$ and $\underline{\bm{\mu}}=((\mu_{1}^{(1)},\ldots,\mu_{d}^{(1)}),$ $\ldots,(\mu_{1}^{(r)},\ldots,\mu_{d}^{(r)}))$ be two $(r,d)$-partitions of $n.$ We say that $\underline{\bm{\lambda}}$ dominates $\underline{\bm{\mu}},$ and we write $\underline{\bm{\lambda}}\unrhd \underline{\bm{\mu}}$ if and only if $$\sum_{i=1}^{k-1}\sum_{j=1}^{d}|\lambda_{j}^{(i)}|+\sum_{j=1}^{l-1}|\lambda_{j}^{(k)}|+\sum_{i=1}^{p}\lambda_{l,i}^{(k)}\geq \sum_{i=1}^{k-1}\sum_{j=1}^{d}|\mu_{j}^{(i)}|+\sum_{j=1}^{l-1}|\mu_{j}^{(k)}|+\sum_{i=1}^{p}\mu_{l,i}^{(k)}$$
for all $k,$ $l$ and $p$ with $1\leq k\leq r,$ $1\leq l\leq d$ and $p\geq 0.$ If $\underline{\bm{\lambda}}\unrhd\underline{\bm{\mu}}$ and $\underline{\bm{\lambda}}\neq \underline{\bm{\mu}},$ we write $\underline{\bm{\lambda}}\rhd \underline{\bm{\mu}}.$
\end{definition}

We extend the partial order above to standard $(r,d)$-tableaux as follows. If $\mathfrak{v}$ is a row standard $(r,d)$-tableau of shape $\underline{\bm{\lambda}}$ and $1\leq k\leq n,$ then the entries $1,2,\ldots,k$ in $\mathfrak{v}$ occupy the diagrams of an $(r,d)$-composition; let $\mathrm{shape}(\mathfrak{v}_{\downarrow k})$ denote this $(r,d)$-composition. Let $\underline{\bm{\lambda}}, \underline{\bm{\mu}}\in \mathcal{C}_{r,n}^{d}.$ Suppose that $\mathfrak{s}$ is a row standard $(r,d)$-tableau of shape $\underline{\bm{\lambda}}$ and that $\mathfrak{t}$ is a row standard $(r,d)$-tableau of shape $\underline{\bm{\mu}}$. We say that $\mathfrak{s}$ dominates $\mathfrak{t},$ and we write $\mathfrak{s}\unrhd\mathfrak{t}$ if $\mathrm{shape}(\mathfrak{s}_{\downarrow k})\unrhd \mathrm{shape}(\mathfrak{t}_{\downarrow k})$ for all $k.$ If $\mathfrak{s}\unrhd\mathfrak{t}$ and $\mathfrak{s}\neq\mathfrak{t},$ then we write $\mathfrak{s}\rhd\mathfrak{t}.$

\subsection{Cellular bases}
Let $\zeta=e^{2\pi i/r}.$ We now fix once and for all a total order on the set of $r$-th roots of unity via setting $\zeta_{k} :=\zeta^{k-1}$ for $1\leq k\leq r.$ Set $S :=\{\zeta_{1},\zeta_{2},\ldots,\zeta_{r}\}.$ Then we define a set partition $A_{\underline{\bm{\lambda}}}\in \mathcal{SP}_{n}$ for any $(r,d)$-partition $\underline{\bm{\lambda}}.$

\begin{definition}\label{definition-1}
Let $\underline{\bm{\lambda}}=((\lambda_{1}^{(1)},\ldots,\lambda_{d}^{(1)}),\ldots,(\lambda_{1}^{(r)},\ldots,\lambda_{d}^{(r)}))\in \mathcal{P}_{r,n}^{d}.$ Suppose that we choose all $1\leq i_{1}< i_{2}<\cdots < i_{p}\leq r$ such that $(\lambda_{1}^{(i_1)},\ldots,\lambda_{d}^{(i_1)}),$ $(\lambda_{1}^{(i_2)},\ldots,\lambda_{d}^{(i_2)}),\ldots,$$(\lambda_{1}^{(i_p)},$\\$\ldots,\lambda_{d}^{(i_p)})$ are nonempty. Define $a_{k} :=\sum_{j=1}^{k}|\bm{\lambda}^{(i_{j})}|$ for $1\leq k\leq p,$ where $|\bm{\lambda}^{(i_{j})}|=\sum_{l=1}^{d}|\lambda_{l}^{(i_{j})}|.$ Then the set partition $A_{\underline{\bm{\lambda}}}$ associated with $\underline{\bm{\lambda}}$ is defined as $$A_{\underline{\bm{\lambda}}} :=\{\{1,\ldots,a_{1}\},\{a_{1}+1,\ldots,a_{2}\},\ldots,\{a_{p-1}+1,\ldots,n\}\},$$ which may be written as $A_{\underline{\bm{\lambda}}}=\{I_{1},I_{2},\ldots,I_{p}\},$ and is referred to the blocks of $A_{\underline{\bm{\lambda}}}$ in the order given above.
\end{definition}

\begin{definition}\label{definition-2}
Let $\underline{\bm{\lambda}}=((\lambda_{1}^{(1)},\ldots,\lambda_{d}^{(1)}),\ldots,(\lambda_{1}^{(r)},\ldots,\lambda_{d}^{(r)}))\in \mathcal{P}_{r,n}^{d},$ and let $a_{k} :=\sum_{j=1}^{k}|\bm{\lambda}^{(i_{j})}|$ $(1\leq k\leq p)$ be defined as above. Then we define $$u_{\underline{\bm{\lambda}}} :=u_{a_{1},i_{1}}u_{a_{2},i_{2}}\cdots u_{a_{p},i_{p}},$$ where $u_{i,k}=\Pi_{l=1;l\neq k}^{r}(t_{i}-\zeta_{l})$ for $1\leq i\leq n$ and $1\leq k\leq r.$
\end{definition}

\begin{definition}\label{definition-3}
Let $\underline{\bm{\lambda}}=((\lambda_{1}^{(1)},\ldots,\lambda_{d}^{(1)}),\ldots,(\lambda_{1}^{(r)},\ldots,\lambda_{d}^{(r)}))\in \mathcal{P}_{r,n}^{d}.$ Associated with $\underline{\bm{\lambda}}$ we can define the following elements $a_{l}^{k}$ and $b_{k}$:$$a_{l}^{k} :=\sum_{m=1}^{l-1}|\lambda_{m}^{(k)}|,~~~~b_{k} :=\sum_{j=1}^{k-1}\sum_{i=1}^{d}|\lambda_{i}^{(j)}|~~~~\mathrm{for}~1\leq k\leq r~\mathrm{and}~1\leq l\leq d.$$ Associated with these elements we can define an element $u_{\mathbf{a}}^{+} :=u_{\mathbf{a}, 1}u_{\mathbf{a}, 2}\cdots u_{\mathbf{a}, r},$ where $$u_{\mathbf{a}, k} :=\prod_{l=1}^{d}\prod_{j=1}^{a_{l}^{k}}(X_{b_{k}+j}-v_{l}).$$
\end{definition}

\begin{definition}\label{definition-4}
Let $\underline{\bm{\lambda}}\in \mathcal{P}_{r,n}^{d}$ and define $u_{\mathbf{a}}^{+}$ as above. We set $U_{\underline{\bm{\lambda}}} :=u_{\underline{\bm{\lambda}}}E_{A_{\underline{\bm{\lambda}}}},$ and define $x_{\underline{\bm{\lambda}}}=\sum_{w\in \mathfrak{S}_{\underline{\bm{\lambda}}}}g_{w}.$ Then we define the element $m_{\underline{\bm{\lambda}}}$ of $Y_{r,n}^{d}$ as follows: $$m_{\underline{\bm{\lambda}}} :=U_{\underline{\bm{\lambda}}}u_{\mathbf{a}}^{+}x_{\underline{\bm{\lambda}}}=
u_{\underline{\bm{\lambda}}}E_{A_{\underline{\bm{\lambda}}}}u_{\mathbf{a}}^{+}x_{\underline{\bm{\lambda}}}.$$
\end{definition}

Let $\ast$ denote the $\mathbb{K}$-linear anti-automorphism of $Y_{r,n}^{d}$, which is determined by $$g_{i}^{\ast}=g_{i},~~~~t_{j}^{\ast}=t_{j},~~~~x_{j}^{\ast}=x_{j}~~~\mathrm{for}~1\leq i\leq n-1~\mathrm{and}~1\leq j\leq n.$$

\begin{definition}\label{definition-5}
Let $\underline{\bm{\lambda}}\in \mathcal{P}_{r,n}^{d},$ and let $\mathfrak{s}$ and $\mathfrak{t}$ be two row standard $(r,d)$-tableaux of shape $\underline{\bm{\lambda}}.$ We then define $m_{\mathfrak{s}\mathfrak{t}}=g_{d(\mathfrak{s})}^{\ast}m_{\underline{\bm{\lambda}}}g_{d(\mathfrak{t})}.$
\end{definition}

For each $\underline{\bm{\mu}}\in \mathcal{P}_{r,n}^{d},$ let $Y_{r,n}^{d, \rhd \underline{\bm{\mu}}}$ be the $\mathbb{K}$-submodule of $Y_{r,n}^{d}$ spanned by $m_{\mathfrak{u}\mathfrak{v}}$ with $\mathfrak{u}, \mathfrak{v}\in \text{Std}(\underline{\bm{\lambda}})$ for various $\underline{\bm{\lambda}}\in \mathcal{P}_{r,n}^{d}$ such that $\underline{\bm{\lambda}}\rhd \underline{\bm{\mu}}.$

\begin{theorem}\label{cellular-basis-the}
The algebra $Y_{r,n}^{d}$ is a free $\mathbb{K}$-module with a cellular basis $$\mathcal{B}_{r,n}^{d}=\{m_{\mathfrak{s}\mathfrak{t}}\:|\:\mathfrak{s}, \mathfrak{t}\in \emph{Std}(\underline{\bm{\lambda}})~for~some~(r,d)\mathrm{-}partition~\underline{\bm{\lambda}}~of~n\},$$ that is, the following properties hold$:$

$\mathrm{(i)}$ The $\mathbb{K}$-linear map determined by $m_{\mathfrak{s}\mathfrak{t}}\mapsto m_{\mathfrak{t}\mathfrak{s}}$ $(m_{\mathfrak{s}\mathfrak{t}}\in \mathcal{B}_{r,n}^{d})$ is an anti-automorphism on $Y_{r,n}^{d}$.

$\mathrm{(ii)}$ For a given $h\in Y_{r,n}^{d}$ and $\mathfrak{t}\in \emph{Std}(\underline{\bm{\mu}}),$ there exist $r_{\mathfrak{v}}\in \mathbb{K}$ such that for all $\mathfrak{s}\in \emph{Std}(\underline{\bm{\mu}})$, we have $$m_{\mathfrak{s}\mathfrak{t}}h\equiv\sum_{\mathfrak{v}\in \emph{Std}(\underline{\bm{\mu}})}r_{\mathfrak{v}}m_{\mathfrak{s}\mathfrak{v}}~~~~\mathrm{mod}~Y_{r,n}^{d, \rhd \underline{\bm{\mu}}},$$
where $r_{\mathfrak{v}}$ may depend on $\mathfrak{v}, \mathfrak{t}$ and $h,$ but not on $\mathfrak{s}.$
\end{theorem}
\begin{proof}
The proof of this result is similar to, but much easier than, the corresponding result for the cyclotomic Yokonuma-Hecke algebras; see [C1, Section 6] for details.
\end{proof}

\subsection{Jucys-Murphy elements}
\begin{definition}\label{definition-6}
Let $\mathfrak{t}$ be an $(r,d)$-tableau of shape $\underline{\bm{\lambda}}$ and suppose that the $(r,d)$-node $\bm{\theta}$ of $\mathfrak{t}$ labelled by $((a, b), k, l)$ is filled in with the element $i$ $(1\leq i\leq n)$. Then we define the content of $i$ as the element $\text{c}_{\mathfrak{t}}(i) :=\text{c}(\bm{\theta})=b-a+v_{l}$ and the $r$-position of $i$ as the element $\text{p}_{\mathfrak{t}}(i) :=\text{p}^{(r)}(\bm{\theta})=k,$ respectively.
\end{definition}

The next proposition claims that the following set
$$\mathcal{L}=\{L_{1},\ldots,L_{2n}\:|\:L_{k}=x_{k}, L_{n+k}=t_{k}~\mathrm{for}~1\leq k\leq n\}$$
is a family of JM-elements for $Y_{r,n}^{d}$ in the abstract sense defined in [Ma, Definition 2.4], which is with respect to the cellular basis in Theorem \ref{cellular-basis-the} .
\begin{proposition}\label{Jucys-Murphy-ele}
Suppose that $\underline{\bm{\lambda}}\in \mathcal{P}_{r,n}^{d}$ and that $\mathfrak{s}, \mathfrak{t}\in \emph{Std}(\underline{\bm{\lambda}}).$ Then we have
\begin{equation}\label{JM-ele-1}
m_{\mathfrak{s}\mathfrak{t}}\hspace{0.3mm}x_k\equiv \emph{c}_{\mathfrak{t}}(k)\hspace{0.3mm}m_{\mathfrak{s}\mathfrak{t}}+\sum_{\substack{\mathfrak{v}\in \emph{Std}(\underline{\bm{\lambda}})\\\mathfrak{v}\rhd \mathfrak{t}}}r_{\mathfrak{v}\mathfrak{t}}\hspace{0.3mm}m_{\mathfrak{s}\mathfrak{v}}~~~~\mathrm{mod}~Y_{r,n}^{d, \rhd \underline{\bm{\lambda}}}
\end{equation}
for some $r_{\mathfrak{v}\mathfrak{t}}\in \mathbb{K}.$ And moreover, we have
\begin{equation}\label{JM-ele-2}
m_{\mathfrak{s}\mathfrak{t}}\hspace{0.3mm}t_k=\zeta_{\emph{p}_{\mathfrak{t}}(k)}m_{\mathfrak{s}\mathfrak{t}}.
\end{equation}
\end{proposition}
\begin{proof}
This proposition can be proved in exactly the same way as in [C1, Proposition 7.3]. We omit the details.
\end{proof}

From Proposition \ref{Jucys-Murphy-ele}, we can now apply the general theory developed in [Ma] to recover the following semi-simplicity criterion of $Y_{r,n}^{d}$ obtained in [C4, Section 4(4.1)]:
\begin{equation}\label{semisimple-condition}
n!\prod_{1\leq i< j\leq d}\prod_{-n< l< n}(l+v_i-v_j)\neq 0.
\end{equation}
In the rest of this paper, we shall work with a split semisimple degenerate cyclotomic Yokonuma-Hecke algebra $Y_{r, n}^{d}$ defined over $\mathbb{K},$ that is, $v_i\in \mathbb{K},$ for $1\leq i\leq d,$ satisfies the condition \eqref{semisimple-condition}. In particular, we can apply all the results in [Ma, Section 3].

\section{Fusion procedure for $Y_{r, n}^{d}$}

In this section, for the split semisimple $Y_{r,n}^{d}$, we first deduce some properties of the idempotents $E_{\mathfrak{t}}$ by applying the general approach developed in [Ma, Section 3]. Then we give the inductive formulae for the idempotents $E_{\mathfrak{t}}$ in terms of the Jucys-Murphy elements. Finally we establish the fusion formulae for $E_{\mathfrak{t}}.$

\subsection{Idempotents of $Y_{r, n}^{d}$}
We first follow the arguments of [Ma, Section 3] to define the primitive idempotents $E_{\mathfrak{t}}$, and use them to construct a ``seminormal" basis of $Y_{r, n}^{d}$. For $1\leq k\leq n,$ we define the following two sets:
\[\mathcal{C}(k) :=\{\text{c}_{\mathfrak{t}}(k)\:|\:\mathfrak{t}\in \text{Std}(\underline{\bm{\lambda}})\text{ for some }\underline{\bm{\lambda}}\in \mathcal{P}_{r,n}^{d}\},\]
and
\[\overline{\mathcal{C}(k)} :=\{\zeta_{\text{p}_{\mathfrak{t}}(k)}\:|\:\mathfrak{t}\in \text{Std}(\underline{\bm{\lambda}})\text{ for some }\underline{\bm{\lambda}}\in \mathcal{P}_{r,n}^{d}\}.\]
\begin{definition}\label{definition-idempotent}
Suppose that $\underline{\bm{\lambda}}\in \mathcal{P}_{r,n}^{d}$ and that $\mathfrak{s}, \mathfrak{t}\in \text{Std}(\underline{\bm{\lambda}}).$

(i) Let
\begin{equation}\label{idempotent-ET}
E_{\mathfrak{t}}=\prod_{k=1}^{n}\bigg(\prod_{\substack{c\in \mathcal{C}(k)\\c\neq \text{c}_{\mathfrak{t}}(k)}}\frac{x_{k}-c}{\text{c}_{\mathfrak{t}}(k)-c}\cdot \prod_{\substack{\bar{c}\in \overline{\mathcal{C}(k)}\\\bar{c}\neq \zeta_{\text{p}_{\mathfrak{t}}(k)}}}\frac{t_{k}-\bar{c}}{\zeta_{\text{p}_{\mathfrak{t}}(k)}-\bar{c}}
\bigg).
\end{equation}

(ii) Let $e_{\mathfrak{s}\mathfrak{t}}^{\underline{\bm{\lambda}}}=E_{\mathfrak{s}}\hspace{0.3mm}m_{\mathfrak{s}\mathfrak{t}}\hspace{0.3mm}E_{\mathfrak{t}}.$
\end{definition}

By Proposition \ref{Jucys-Murphy-ele} and applying the general theory developed in [Ma, Section 3], we can get the following results.

\begin{proposition}\label{idempotent-property}
$(\emph{i})$ The set $\{e_{\mathfrak{s}\mathfrak{t}}^{\underline{\bm{\lambda}}}\:|\:\mathfrak{s}, \mathfrak{t}\in \emph{Std}(\underline{\bm{\lambda}})\text{ for some }\underline{\bm{\lambda}}\in \mathcal{P}_{r,n}^{d}\}$ is a $\mathbb{K}$-basis of $Y_{r, n}^{d}.$

$(\emph{ii})$ For $\underline{\bm{\lambda}}, \underline{\bm{\mu}}\in \mathcal{P}_{r,n}^{d}$ and $\mathfrak{s}, \mathfrak{t}\in \emph{Std}(\underline{\bm{\lambda}}),$ $\mathfrak{u}, \mathfrak{v}\in \emph{Std}(\underline{\bm{\mu}}),$ we have
\begin{equation}\label{result-1}
e_{\mathfrak{s}\mathfrak{t}}^{\underline{\bm{\lambda}}}\hspace{0.3mm}x_k=\emph{c}_{\mathfrak{t}}(k)e_{\mathfrak{s}\mathfrak{t}}^{\underline{\bm{\lambda}}},\hspace{1cm}
e_{\mathfrak{s}\mathfrak{t}}^{\underline{\bm{\lambda}}}\hspace{0.3mm}t_k=\zeta_{\emph{p}_{\mathfrak{t}}(k)}e_{\mathfrak{s}\mathfrak{t}}^{\underline{\bm{\lambda}}},\hspace{1cm} e_{\mathfrak{s}\mathfrak{t}}^{\underline{\bm{\lambda}}}\hspace{0.3mm}E_{\mathfrak{u}}=\delta_{\mathfrak{t}, \mathfrak{u}}e_{\mathfrak{s}\mathfrak{u}}^{\underline{\bm{\lambda}}},
\end{equation}
and moreover, there exists a scalar $0\neq \gamma_{\mathfrak{t}}\in \mathbb{K}$ such that
\begin{align}\label{result-2}
e_{\mathfrak{s}\mathfrak{t}}^{\underline{\bm{\lambda}}}\hspace{0.3mm}e_{\mathfrak{u}\mathfrak{v}}^{\underline{\bm{\mu}}}=
\begin{cases}
\gamma_{\mathfrak{t}}\hspace{0.3mm}e_{\mathfrak{s}\mathfrak{v}}^{\underline{\bm{\lambda}}} & \text{if } \underline{\bm{\lambda}}=\underline{\bm{\mu}}\text{ and }\mathfrak{t}=\mathfrak{u};
\\
\hspace{0.25cm}0 & \text{otherwise.}
\end{cases}
\end{align}
In particular, $\gamma_{\mathfrak{t}}$ depends only on $\mathfrak{t}$ and the set $\{e_{\mathfrak{s}\mathfrak{t}}^{\underline{\bm{\lambda}}}\:|\:\mathfrak{s}, \mathfrak{t}\in \emph{Std}(\underline{\bm{\lambda}})\text{ and }\underline{\bm{\lambda}}\in \mathcal{P}_{r,n}^{d}\}$ is a cellular basis of $Y_{r, n}^{d}.$

$(\emph{iii})$ For $\underline{\bm{\lambda}}\in \mathcal{P}_{r,n}^{d}$ and $\mathfrak{t}\in \emph{Std}(\underline{\bm{\lambda}}),$ we have $E_{\mathfrak{t}}=\frac{1}{\gamma_{\mathfrak{t}}}\hspace{0.3mm}e_{\mathfrak{t}\mathfrak{t}}^{\underline{\bm{\lambda}}}.$ Moreover, these elements $\{E_{\mathfrak{t}}\:|\:\mathfrak{t}\in \emph{Std}(\underline{\bm{\lambda}})\text{ for some }\underline{\bm{\lambda}}\in \mathcal{P}_{r,n}^{d}\}\}$ give a complete set of pairwise orthogonal primitive idempotents for $Y_{r, n}^{d}.$

$(\emph{iv})$ For $\underline{\bm{\lambda}}\in \mathcal{P}_{r,n}^{d}$ and $\mathfrak{t}\in \emph{Std}(\underline{\bm{\lambda}}),$ we have
\begin{equation}\label{result-3}
E_{\mathfrak{t}}\hspace{0.3mm}x_k=x_k\hspace{0.3mm}E_{\mathfrak{t}}=\emph{c}_{\mathfrak{t}}(k)E_{\mathfrak{t}},\hspace{1cm}
E_{\mathfrak{t}}\hspace{0.3mm}t_k=t_k\hspace{0.3mm}E_{\mathfrak{t}}=\zeta_{\emph{p}_{\mathfrak{t}}(k)}E_{\mathfrak{t}}.
\end{equation}

$(\emph{v})$ The Jucys-Murphy elements $x_1,\ldots,x_n,$ $t_1,\ldots,t_n$ generate a maximal commutative subalgebra of $Y_{r, n}^{d}.$
\end{proposition}

\subsection{Inductive formulae of $E_{\mathfrak{t}}$}
We first define two rational functions, and then introduce the inductive formulae of the primitive idempotents $E_{\mathfrak{t}}.$

For $\lambda, \mu$ two partitions and a node $\theta=(x,y)\in [\lambda],$ we define the hook length $h_{\lambda}(\theta)$ and the generalized hook length $h_{\lambda}^{\mu}(\theta)$ of $\theta$ with respect to $(\lambda, \mu)$, respectively:
\begin{equation}\label{hook-length1}
h_{\lambda}(\theta) :=\lambda_{x}+\lambda'_{y}-x-y+1,
\end{equation}
and
\begin{equation}\label{hook-length2}
h_{\lambda}^{\mu}(\theta) :=\lambda_{x}+\mu'_{y}-x-y+1.
\end{equation}

Assume that $\underline{\bm{\lambda}}=((\lambda_{1}^{(1)},\ldots,\lambda_{d}^{(1)}),\ldots,(\lambda_{1}^{(r)},\ldots,\lambda_{d}^{(r)}))$ is an $(r,d)$-partition and $\bm{\theta}=(\theta, k,l)=((x,y),k,l)$ is an $(r,d)$-node of $[\underline{\bm{\lambda}}].$ The hook length $h_{\underline{\bm{\lambda}}}(\bm{\theta})$ of $\bm{\theta}$ in $\underline{\bm{\lambda}}$ is defined to be the hook length of the node $\theta$ in the partition of $\underline{\bm{\lambda}}$ with position $(k,l)$, that is,
\begin{equation}\label{hook-length3}
h_{\underline{\bm{\lambda}}}(\bm{\theta}) :=h_{\lambda^{(k)}_{l}}(\theta)=\lambda^{(k)}_{l,x}+\lambda^{(k)'}_{l,y}-x-y+1.
\end{equation}
Let $\mu$ be another partition. We define the generalized hook length $h_{\underline{\bm{\lambda}}}^{\mu}(\bm{\theta})$ of $\bm{\theta}$ with respect to $(\underline{\bm{\lambda}}, \mu)$ to be the generalized hook length of $\theta$ with respect to $(\lambda^{(k)}_{l}, \mu),$ that is,
\begin{equation}\label{hook-length4}
h_{\underline{\bm{\lambda}}}^{\mu}(\bm{\theta}) :=h_{\lambda^{(k)}_{l}}^{\mu}(\theta)=\lambda^{(k)}_{l,x}+\mu'_{y}-x-y+1.
\end{equation}

Recall that $S=\{\zeta_1,\ldots,\zeta_r\}$ is the set of all $r$-th roots of unity. For an $(r,d)$-partition $\underline{\bm{\lambda}}=((\lambda_{1}^{(1)},\ldots,\lambda_{d}^{(1)}),\ldots,(\lambda_{1}^{(r)},\ldots,\lambda_{d}^{(r)}))$, we define
\begin{equation}\label{f-lambda-t}
\text{F}_{\underline{\bm{\lambda}}}^{T} :=\prod_{\bm{\theta}\in \underline{\bm{\lambda}}} \Big(\prod_{\substack{\xi\in S\\\xi\neq \zeta_{\text{p}^{(r)}(\bm{\theta})}}}(\zeta_{\text{p}^{(r)}(\bm{\theta})}-\xi)\Big).
\end{equation}
and
\begin{equation}\label{f-lambda}
\text{F}_{\underline{\bm{\lambda}}} :=\prod_{\bm{\theta}\in \underline{\bm{\lambda}}}
\bigg(h_{\underline{\bm{\lambda}}}(\bm{\theta})\prod_{\substack{1\leq k\leq d\\k\neq \text{p}^{(d)}(\bm{\theta})}}\Big(h_{\underline{\bm{\lambda}}}^{\lambda^{(\text{p}^{(r)}(\bm{\theta}))}_{k}}
(\bm{\theta})+v_{\text{p}^{(d)}(\bm{\theta})}-v_{k}\Big)\bigg).
\end{equation}

\begin{remark}\label{remark1}
When $r=1$ and $d=m,$ the degenerate cyclotomic Yokonuma-Hecke algebra $Y_{r,n}^{d}$ is exactly the cyclotomic Hecke algebra $H_{n}^{m},$ and $\text{F}_{\underline{\bm{\lambda}}}$ is precisely the Schur element of $H_{n}^{m}$ which has been explicitly calculated in [Z, Theorem 4.2].
\end{remark}

Let $\underline{\bm{\lambda}}$ be an $(r,d)$-partition of $n$ and $\mathfrak{t}$ a standard $(r,d)$-tableau of shape $\underline{\bm{\lambda}}$. In order to simplify notation, we set $\text{c}_i :=\text{c}_{\mathfrak{t}}(i)$ and $\text{p}_i :=\text{p}_{\mathfrak{t}}(i)$ for $i=1,\ldots,n.$ We then define
\begin{equation}\label{f-TT}
\text{F}_{\mathfrak{t}}^{T}(\underline{v}) :=\prod_{\substack{\xi\in S\\\xi\neq \zeta_{\text{p}_{n}}}}\frac{1}{\underline{v}-\xi},
\end{equation}
and
\begin{equation}\label{f-T}
\text{F}_{\mathfrak{t}}(u) :=\frac{u-\text{c}_n}{(u-v_1)\cdots (u-v_d)}\prod_{i=1}^{n-1}\frac{(u-\text{c}_i)^{2}}{(u-\text{c}_i)^{2}-\delta_{\text{p}_i, \text{p}_n}},
\end{equation}
where $\delta_{\text{p}_i, \text{p}_n}$ is the Kronecker delta.

Denote by $\bm{\theta}$ the $(r,d)$-node of $\mathfrak{t}$ containing the number $n.$ Since $\mathfrak{t}$ is standard, the $(r,d)$-node $\bm{\theta}$ is removable from $\underline{\bm{\lambda}}.$ Let $\mathfrak{u}$ be the standard $(r,d)$-tableau obtained from $\mathfrak{t}$ by removing $\bm{\theta}$ and let $\underline{\bm{\mu}}$ be the shape of $\mathfrak{u}$.

For each $\xi\in S,$ we have
\begin{equation}
\label{young-module-5-19}\prod_{\xi\neq \alpha\in S}(\xi-\alpha)=r\xi^{-1}.
\end{equation}
Thus, $\text{F}_{\mathfrak{t}}^{T}(\underline{v})$ is non-singular at $\underline{v}=\zeta_{\text{p}_{n}}.$ Moreover, by \eqref{f-lambda-t} we get
\begin{equation}\label{f-T-relation}
\text{F}_{\mathfrak{t}}^{T}(\underline{v})\Big|_{\underline{v}=\zeta_{\text{p}_{n}}}=\frac{\zeta_{\text{p}_{n}}}{r}=
(\text{F}_{\underline{\bm{\lambda}}}^{T})^{-1}\text{F}_{\underline{\bm{\mu}}}^{T}.
\end{equation}

The following proposition can be proved in exactly the same way as in [OgPA1, Propositions 3.4 and ZL, Lemma 4.3].
\begin{proposition}
The rational function $\emph{F}_{\mathfrak{t}}(u)$ is non-singular at $u=\emph{c}_n$, and moreover, we have
\begin{equation}\label{f-T-u}
\emph{F}_{\mathfrak{t}}(u)\Big|_{u=\emph{c}_n}=\emph{F}_{\underline{\bm{\lambda}}}^{-1}\emph{F}_{\underline{\bm{\mu}}}.
\end{equation}
\end{proposition}

Denote by $\mathcal{E}_{+}(\underline{\bm{\mu}})$ the set of $(r,d)$-nodes addable to $\underline{\bm{\mu}}.$ By \eqref{idempotent-ET}, we have the following inductive formula for $E_{\mathfrak{t}}$:
\begin{equation}\label{inductive-formula}
E_{\mathfrak{t}}=E_{\mathfrak{u}}\prod_{\substack{\bm{\theta'}\in \mathcal{E}_{+}(\underline{\bm{\mu}})\\\text{c}(\bm{\theta'})\neq\text{c}(\bm{\theta})}}
\frac{x_{n}-\text{c}(\bm{\theta'})}{\text{c}(\bm{\theta})-\text{c}(\bm{\theta'})}\prod_{\substack{\bm{\theta'}\in \mathcal{E}_{+}(\underline{\bm{\mu}})\\\text{p}^{(r)}(\bm{\theta'})\neq\text{p}^{(r)}(\bm{\theta})}}
\frac{t_{n}-\zeta_{\text{p}^{(r)}(\bm{\theta'})}}{\zeta_{\text{p}^{(r)}(\bm{\theta})}-\zeta_{\text{p}^{(r)}(\bm{\theta'})}}
\end{equation}
with $E_{\mathfrak{t}_{0}}=1$ for the unique standard $(r,d)$-tableau $\mathfrak{t}_{0}$ of size $0.$

\subsection{Fusion formulae of $E_{\mathfrak{t}}$}
We shall present a fusion formula for the idempotent $E_{\mathfrak{t}}.$ We first define the following rational functions in variables $a,b$ with values in $Y_{r, n}^{d}$:
\begin{equation}\label{Baxter-element}
g_{i}(a,b) :=g_{i}+\frac{e_{i}}{a-b}\quad\mbox{for}~i=1,\ldots,n-1.
\end{equation}
The following lemma is proved in [PA1, Proposition 2].
\begin{lemma}\label{Bax-elements}
The rational functions $g_{i}(a,b)$ satisfy the following relations$:$
\begin{equation}\label{Baxter-element1}
g_{i}(a,b)g_{i+1}(a,c)g_{i}(b,c)=g_{i+1}(b,c)g_{i}(a,c)g_{i+1}(a,b)\quad\mbox{for}~i=1,\ldots,n-1,
\end{equation}
\begin{equation}\label{Baxter-element2}
\hspace*{20pt}g_{i}(a,b)g_{i}(b,a)=1-\frac{e_{i}}{(a-b)^{2}}\qquad\mbox{for}~i=1,\ldots,n-1.
\end{equation}
\end{lemma}

Set $\phi_{1}(u) :=\frac{(u-v_1)\cdots (u-v_d)}{u-x_1}.$ For $k=2,\ldots,n$, we set
\begin{align}\label{phi-function}
\phi_k(u_1,\ldots,u_{k-1},u) :&=g_{k-1}(u,u_{k-1})\phi_{k-1}(u_1,\ldots,u_{k-2},u)g_{k-1}\notag\\
&=g_{k-1}(u,u_{k-1})g_{k-2}(u,u_{k-2})\cdots g_{1}(u,u_{1})\phi_{1}(u)\cdot g_{1}\cdots g_{k-1}.
\end{align}

Recall that $S=\{\zeta_1,\ldots,\zeta_r\}$ is the set of all $r$-th roots of unity. We set
\begin{equation}\label{gamma-function}
\Gamma(\underline{v_1},\ldots,\underline{v_n}) :=\prod_{i=1}^{n}\Big(\frac{\Pi_{\xi\in S}(\underline{v_i}-\xi)}{\underline{v_i}-t_i}\Big).
\end{equation}

We then define the following rational function:
\begin{align}\label{Phi-function}
\Phi(u_1,\ldots,u_n,\underline{v_1},\ldots,\underline{v_n}) :=&\phi_n(u_1,\ldots,u_{n})\phi_{n-1}(u_1,\ldots,u_{n-1})\notag\\
&\cdots \phi_1(u_1)\Gamma(\underline{v_1},\ldots,\underline{v_n}).
\end{align}

Now we can state the main result of this section.

\begin{theorem}\label{main-theorem}
The idempotent $E_{\mathfrak{t}}$ of $Y_{r,n}^{d}$ corresponding to the standard $(r,d)$-tableau $\mathfrak{t}$ can be obtained by the following consecutive evaluations$:$
\begin{equation}\label{idempotents}
E_{\mathfrak{t}}=\frac{1}{\emph{F}_{\underline{\bm{\lambda}}}^{T}\emph{F}_{\underline{\bm{\lambda}}}}
\Phi(u_1,\ldots,u_{n},\underline{v_1},\ldots,\underline{v_n})
\Big|_{\underline{v_1}=\zeta_{\emph{p}_{1}}}\cdots\Big|_{\underline{v_n}=\zeta_{\emph{p}_{n}}}\Big|_{u_{1}=\emph{c}_1}\cdots\Big|_{u_{n}=\emph{c}_{n}}.
\end{equation}
\end{theorem}

We first prove two necessary lemmas. We define the following element:
\begin{equation}\label{e-upn}
E_{\mathfrak{u}, \text{p}_{n}} :=\frac{\underline{v}-\zeta_{\text{p}_n}}{\underline{v}-t_{n}}E_{\mathfrak{u}}\Big|_{\underline{v}=\zeta_{\text{p}_{n}}}.
\end{equation}
Then the element $E_{\mathfrak{u}, \text{p}_{n}}$ is an idempotent which is equal to the sum of the idempotents $E_{\mathfrak{v}}$, where $\mathfrak{v}$ runs through the set of standard $(r,d)$-tableaux obtained from $\mathfrak{u}$ by adding an $(r,d)$-node $\bm{\theta}$ containing the integer $n$ such that $\text{p}^{(r)}(\bm{\theta})=\text{p}_{n}.$

\begin{lemma}
Assume that $n\geq 1.$ We have
\begin{align}\label{F-PhiEu}
\emph{F}_{\mathfrak{t}}(u)\phi_{n}(\emph{c}_1,\ldots,\emph{c}_{n-1},u)E_{\mathfrak{u}, \emph{p}_{n}}=\frac{u-\emph{c}_n}{u-x_{n}}E_{\mathfrak{u}, \emph{p}_{n}}.
\end{align}
\end{lemma}
\begin{proof}
The proof of this result is similar to, but easier than, the corresponding result for the cyclotomic Yokonuma-Hecke algebras; see [C6] for details.
\end{proof}

For $k=1,\ldots,n,$ we define
\begin{equation}\label{tilde-phi}
\widetilde{\phi}_{k}(u_1,\ldots,u_{k-1},u,\underline{v}) :=\phi_k(u_1,\ldots,u_{k-1},u)\cdot\Big(\frac{\Pi_{\xi\in S}(\underline{v}-\xi)}{\underline{v}-t_k}\Big).
\end{equation}

\begin{lemma}
Assume that $n\geq 1.$ We have
\begin{align}\label{FF-Phi}
\emph{F}_{\mathfrak{t}}^{T}(\underline{v})\emph{F}_{\mathfrak{t}}(u)\widetilde{\phi}_{n}(\emph{c}_1,\ldots,\emph{c}_{n-1},u,\underline{v})
E_{\mathfrak{u}}\Big|_{\underline{v}=\zeta_{\emph{p}_{n}}}=\frac{u-\emph{c}_n}{u-x_{n}}
\frac{\underline{v}-\zeta_{\emph{p}_n}}{\underline{v}-t_{n}}E_{\mathfrak{u}}\Big|_{\underline{v}=\zeta_{\emph{p}_{n}}}.
\end{align}
\end{lemma}
\begin{proof}
By \eqref{f-TT}, we have
\begin{align}\label{Ftt-equality}
\emph{F}_{\mathfrak{t}}^{T}(\underline{v})\cdot\Big(\frac{\Pi_{\xi\in S}(\underline{v}-\xi)}{\underline{v}-t_n}\Big)=\frac{\underline{v}-\zeta_{\text{p}_n}}{\underline{v}-t_{n}}.
\end{align}
By \eqref{e-upn}, \eqref{tilde-phi} and \eqref{Ftt-equality}, we see that \eqref{FF-Phi} is a direct consequence of \eqref{F-PhiEu}.
\end{proof}

{\it Proof of Theorem \ref{main-theorem}} Since $g_i$ commutes with $t_k$ if $i< k-1$, we can rewrite the rational function $\Phi(u_1,\ldots,u_n,v_1,\ldots,v_n)$ as follows:
\begin{align}\label{Phitilde}
\Phi(&u_1,\ldots,u_n,\underline{v_1},\ldots,\underline{v_n})\notag\\
&=\widetilde{\phi}_{n}(u_1,\ldots,u_{n},\underline{v_n})\widetilde{\phi}_{n-1}(u_1,\ldots,u_{n-1},\underline{v_{n-1}})\cdots \widetilde{\phi}_{1}(u_1,\underline{v_1}).
\end{align}

We shall prove this theorem by induction on $n.$ For $n=0,$ the situation is trivial.

For $n> 0,$ by \eqref{Phitilde} and the induction hypothesis we can rewrite the right-hand side of \eqref{idempotents} as follows:
\begin{equation}\label{ffff-phieu}
(\text{F}_{\underline{\bm{\lambda}}}^{T}\text{F}_{\underline{\bm{\lambda}}})^{-1}\text{F}_{\underline{\bm{\mu}}}^{T}\text{F}_{\underline{\bm{\mu}}}
\widetilde{\phi}_{n}(\text{c}_1,\ldots,\text{c}_{n-1},u_n,\underline{v_n})E_{\mathfrak{u}}\Big|_{\underline{v_n}=\zeta_{\text{p}_{n}}}\Big|_{u_n=\text{c}_{n}}.
\end{equation}
By \eqref{FF-Phi} we can rewrite the expression \eqref{ffff-phieu} as
\begin{equation}\label{ffff-phieu1}
(\text{F}_{\underline{\bm{\lambda}}}^{T}\text{F}_{\underline{\bm{\lambda}}})^{-1}\text{F}_{\underline{\bm{\mu}}}^{T}\text{F}_{\underline{\bm{\mu}}}
(\text{F}_{\mathfrak{t}}^{T}(\underline{v_n})\text{F}_{\mathfrak{t}}(u_n))^{-1}
\frac{u_n-\text{c}_n}{u_n-x_{n}}\frac{\underline{v_n}-\zeta_{\text{p}_n}}{\underline{v_n}-t_{n}}
E_{\mathfrak{u}}\Big|_{\underline{v_n}=\zeta_{\text{p}_{n}}}\Big|_{u_n=\text{c}_{n}}.
\end{equation}

Let $\{\mathfrak{t}_1,\ldots,\mathfrak{t}_k\}$ be the set of pairwise different standard $(r,d)$-tableaux obtained from $\mathfrak{u}$ by adding an $(r,d)$-node containing the integer $n.$ Note that $\mathfrak{t}\in \{\mathfrak{t}_1,\ldots,\mathfrak{t}_k\}.$ Moreover, by branching properties, we have
\begin{equation}\label{sum-formula}
E_{\mathfrak{u}}=\sum_{i=1}^{k}E_{\mathfrak{t}_{i}}.
\end{equation}
If we consider the following rational function in $u$ and $\underline{v}$:
\begin{equation}\label{rational-function}
\frac{u-\text{c}_n}{u-x_{n}}\frac{\underline{v}-\zeta_{\text{p}_n}}{\underline{v}-t_{n}}E_{\mathfrak{u}},
\end{equation}
the formulae \eqref{result-3} imply that \eqref{rational-function} is non-singular at $u=\text{c}_n$ and $\underline{v}=\zeta_{\text{p}_n}$; and moreover, by replacing $E_{\mathfrak{u}}$ with the right-hand side of \eqref{sum-formula}, we get
\begin{equation}\label{sum-function}
\frac{u-\text{c}_n}{u-x_{n}}\frac{\underline{v}-\zeta_{\text{p}_n}}{\underline{v}-t_{n}}E_{\mathfrak{u}}
\Big|_{\underline{v}=\zeta_{\text{p}_{n}}}\Big|_{u=\text{c}_n}=E_{\mathfrak{t}}.
\end{equation}

By \eqref{f-T-relation} and \eqref{f-T-u}, together with \eqref{ffff-phieu1} and \eqref{sum-function}, we see that the right-hand side of \eqref{idempotents} is equal to $E_{\mathfrak{t}}.$  $\hfill{} \Box$

\begin{example}
Consider the situation that $r=d=2, n=4$ and the $(2,2)$-partition $\underline{\bm{\lambda}}=(((2),(0)),((1)),(1))$ of $4$. We shall consider the following standard $(2,2)$-tableau of shape $\underline{\bm{\lambda}}$:
\[\mathfrak{t}=\left(\left(\hspace{-0.1cm}\begin{array}{l}\fbox{1}\fbox{3}\\[-0.05em] \end{array}
\, ,\, \varnothing\right)
\, ,\,\left(\hspace{-0.1cm}\begin{array}{l}\fbox{2}\\[-0.05em] \end{array}\,,\,\begin{array}{l}\fbox{4}\\[-0.05em] \end{array}\hspace{-0.1cm}\right)\right).\]

Theorem \ref{main-theorem} implies that the idempotent $E_{\mathfrak{t}}$ can be expressed as
\begin{align*}
E_{\mathfrak{t}}=&\frac{\zeta_{1}^{2}\zeta_{2}^{2}}{32(v_1-v_2)(v_2-v_1+1)(v_1-v_2+1)^{2}}\\
&\times g_{3}(v_2,v_1+1)g_{2}(v_2,v_1)g_{1}(v_2,v_1)\phi_{1}(v_2)g_{1}g_{2}g_{3}\\
&\times g_{2}(v_1+1, v_1)g_{1}(v_1+1, v_1)\phi_{1}(v_{1}+1)g_{1}g_{2}\times g_{1}(v_1,v_1)\phi_{1}(v_1)g_{1}\\
&\times\phi_{1}(v_1)(\zeta_1+t_1)(\zeta_2+t_2)(\zeta_1+t_3)(\zeta_2+t_4).
\end{align*}

\end{example}



\end{document}